\documentclass[a4paper,12pt]{article}
\usepackage[utf8]{inputenc}

\usepackage{amsmath}
\usepackage{amssymb}
\usepackage{amsthm}
\usepackage{mathtools}

\usepackage{cite}

\usepackage{booktabs}

\usepackage{esint}

\usepackage{mathrsfs}
\usepackage{bbm}

\usepackage{accents}

\usepackage{xcolor}
\usepackage{graphicx}        

\newcommand{\pt}{\partial}
\newcommand{\R}{\mathbb{R}}

\usepackage{listings}
\title{Time stepping adaptation for subdiffusion problems with non-smooth right-hand sides}
\author{Sebastian Franz\footnote{
          Institute of Scientific Computing, Technische Universit\"at Dresden, Germany.
          \mbox{e-mail}: sebastian.franz@tu-dresden.de},\and
        Natalia Kopteva\footnote{
          Department of Mathematics and Statistics, University of Limerick, Ireland.
          \mbox{e-mail}: natalia.kopteva@ul.ie}
       }
\date{\today}

\begin{document}

\maketitle

\begin{abstract}
  We consider a time-fractional subdiffusion equation with a Caputo derivative in time,
  a general second-order elliptic spatial operator, 
   and a right-hand side that is non-smooth in time.
  The presence of the latter may lead to locking problems in our time stepping
  procedure recently introduced in \cite{Kopteva22,FrK23}. Hence,
  a generalized version of the residual barrier is proposed to rectify the issue.
  We also consider related alternatives to this generalized algorithm, and, furthermore, show that this new residual barrier
  may be useful in the case of a negative reaction coefficient.
\end{abstract}

\section{Introduction}
  We consider time-fractional parabolic equations of the form
  \begin{equation}\label{eq:problem}
    (\pt_t^\alpha +L)u = f\quad\text{in }(0,T)\times\Omega,
  \end{equation}
  posed in the spatial domain $\Omega\subset\R^d$, $d\in\{1,\,2,\,3\}$,
  subject to the initial condition $u(\cdot,0)=u_0(\cdot)$ and homogeneous boundary conditions on $\partial\Omega$.
  Here $L$ is  a general second-order elliptic operator with variable coefficients, and $\pt_t^\alpha$ is the Caputo fractional
  derivative in time, defined for $\alpha\in(0,1)$ and $t>0$, see also \cite{Diet10}, by
  \begin{equation}\label{eq:CaputoEquiv}
    \pt_t^{\alpha} u(\cdot,t) := \frac1{\Gamma(1-\alpha)} \int_{0}^t(t-s)^{-\alpha}\, \pt_s u(\cdot, s)\, ds,
  \end{equation}
  where $\Gamma(\cdot)$ is the Gamma function, and $\pt_s$ denotes the 
  partial derivative in~$s$.

  If the right-hand side $f$ is smooth, a typical solution to \eqref{eq:problem} exhibits an initial singularity of type $t^{\alpha}$. Hence, one efficient way of obtaining
  reliable numerical approximations for such problems is to employ suitable non-uniform temporal meshes, which may be constructed a priori (using appropriate mesh grading
  \cite{SORG17,KopMC19}) or a posteriori, based on the a-posteriori error estimation and adaptive time stepping \cite{Kopteva22,FrK23,Kopt_Stynes_apost22}.
  Importantly, the latter methodology relies on the theoretical a-posteriori error estimation proposed in \cite{Kopteva22}; hence such adaptive algorithms yield reliable computed solutions for arbitrarily large times.

  The stable and efficient implementation of such time stepping algorithms
was specifically addressed \cite{FrK23}
        in the context of higher-order methods, including continuous collocation methods of arbitrary order.
    In contrast to a-priori-chosen meshes, the adaptive algorithm
    was shown to be capable of
capturing both initial singularities and local shocks/peaks in the solution.
For example (see \cite[Example 6.2]{FrK23}), for
  $
    f(x,t)=(1-t)\cdot\sin((x\pi)^2)+t\cdot\exp(-100\cdot(2t-1)^2)
  $,
  with a localised Gaussian pulse at $t=0.5$,
   our time stepping algorithm produced a suitable mesh, resolving both the initial singularity and the local phenomena,
  with the error guaranteed to be below any desired tolerance $TOL$.

  However, if the
  right-hand side becomes less smooth and exhibits discontinuities in time, we have discovered that the algorithm in \cite{FrK23} locks on the approach to such points.
  The purpose of this paper is to rectify the issue by employing the same general methodology, but with
 an appropriately-generalized version of the residual barrier, as described in section~\ref{sec:All-in}.
 This generalized barrier takes into account the location of the singularities in the right-hand side,
 while if such problematic points in time are unknown a priori,
 in section~\ref{sec:heuristics} we discuss their  automatic computation.
 Furthermore,
 two  related alternatives for the generalized time stepping algorithm are presented  in section~\ref{sec:Split}.
  Section~\ref{sec_lambda} demonstrates that exactly the same generalized residual barrier may rectify locking issues in the case of negative reaction coefficient.
  In the final section~\ref{sec:stab} we include some considerations and implementation advices on possible stability problems.

  \section{Generalized Residual Barrier for Interior Singularities}\label{sec:All-in}
  The mesh adaptation algorithm in \cite{FrK23} is based on the a-posteriori error estimation of \cite{Kopteva22}, that can be summarised
  as
  \[
    \|Res(\cdot, t)\|<TOL\cdot \mathcal{R}(t)\;\;\forall\, t>0\quad\Rightarrow\quad \|u-u_h\|\leq TOL\cdot \mathcal{E}(t)\;\;\forall\, t>0,
  \]
  where
  $u_h$ is the numerical approximation with the residual
  $Res:=f-(\pt_t^\alpha+L)u_h$, $TOL$ is the desired tolerance,
  and
  $\|\cdot\|$ a suitable norm (the $L_2(\Omega)$ and $L_\infty(\Omega)$ norms were considered).
  The algorithm hinges on an appropriate choice of the desired error barrier
  $\mathcal{E}(t)$ and the corresponding residual barrier
  $\mathcal{R}(t)$, which should remain positive $\forall\,t>0$, and are related by a simple equation
  $(\pt_t^\alpha+\lambda)\mathcal{E}(t)=\mathcal{R}(t)$ $\forall\,t>0$, where the constant $\lambda$ depends on the spatial operator $L$; see \cite{Kopteva22,FrK23} for details.

In this paper, we shall restrict consideration to the $L_\infty(\Omega)$ norm,
 so $\lambda:=\inf  L[1]$,
 and
$\mathcal{E}$ of type
$\mathcal{E}(t)=1$ for $t>0$ (with $\mathcal{E}(0)=0$), which corresponds \cite{Kopteva22,FrK23}
to
\begin{equation}\label{mathcal_R}
    \mathcal{R}(t) = \lambda+\Gamma(1-\alpha)^{-1}t^{-\alpha}.
\end{equation}
It will be convenient to adapt the convention that $\mathcal{R}(t):=0$
 for $t\le 0$.

  We are interested in the right-hand sides $f$ of type
\begin{equation}\label{f_sum}
    f(x,t)=\sum_{k=0}^KH(t-s_k)\,f_k(x,t),
\end{equation}
  where
  $0=s_0< s_1<s_2<\cdots<s_K<T$, $H(t)$ is the Heaviside step function,
  while each $f_k$ is at least continuous. Thus,
   $f$ may have up to $K$ jumps in $(0,T)$.
   As already mentioned, we have discovered that
   if $f$ exhibits discontinuities in time, the algorithm in \cite{FrK23} locks on the approach to such points.
  To rectify this, we propose the following
  generalized  barrier function $\mathcal{B}(t)$
  for the residual $Res$.

  If we assume that each
  $f_k(x,t)$ is of type $(t-s_k)^{\gamma_k}$ for $t>s_k$, with some $\gamma_k\ge 0$,
  one easily concludes that $u$ is expected to have a singularity of type
  $(t-s_k)^{\gamma_k+\alpha}$ as $t\rightarrow s_k^+$.
  While if all $\gamma_k=0$, then $u$ will exhibit a singularity of type $(t-s_k)^{\alpha}$
  (similar to a typical initial singularity!) as $t\rightarrow s_k^+$.
  This observation implies that the algorithm should be modified so that such multiple interior singularities
are effectively treated in exactly the same way as we have treated the initial singularity
in \cite{Kopteva22,FrK23,Kopt_Stynes_apost22}.
This goal is easily attained by a simple generalization of
  the residual barrier \eqref{mathcal_R} to 
  \begin{gather}\label{eq:newBound}
    \mathcal{B}(t)=\sum_{k=0}^K w\,_k H(t-s_k)\,\mathcal{R}(t-s_k),
  \end{gather}
 where
 $\mathcal{R}(t)$ is defined in \eqref{mathcal_R}, while
 $w_k>0$ are adjustable weights with $\sum_{k=0}^K w_k = W$
 (e.g., $w_k:=1$ $\forall k$ yields $W=K+1$, while $w_k:=2^{-k}$ yields $W\leq2$).

  As a consequence, \cite[Corollary~2.3]{Kopteva22} yields the new a posteriori error bound
  \[
    \|Res(t)\|\leq TOL\cdot \mathcal{B}(t)\;\;\forall\,t>0
    \;\;\Rightarrow\;\;
    \|u-u_h\|
    \le TOL\cdot \sum_{k=0}^K w_k\, H(t-s_k)
     \;\;\forall\,t>0,
  \]
  which immediately implies $\|u-u_h\|\leq TOL\cdot W$.

  The resulting time stepping algorithm is presented as Algorithm~\ref{alg:mesh} below.
  \begin{lstlisting}[float,caption=Adaptive Algorithm,  emph={T_cmp,S}, basicstyle=\footnotesize, escapechar=@, label=alg:mesh]
k := 1;  uh(1) := u0;  mesh(1:2) := [0,tau_init]; % init
while mesh(k)<T
  k    := k+1;
  flag := 0;
  while mesh(k)-mesh(k-1) > tau_min
    uh(k)      := computeSolution(mesh(1:k),uh);
    Res        := computeResidual(uh,mesh);
    ResBarrier := computeResidualBarrier(mesh,S);
    T_cmp      := min(S(S>mesh(k)),T);    % next problem point
    if all(Res<TOL*ResBarrier)            % residual small enough
      if mesh(k)>=T_cmp                   % accept
        break                             % finish or next step
      else                                % ok
        if flag = 2                       % from larger step
          mesh(k+1) := min(mesh(k)+(mesh(k)-mesh(k-1)),T_cmp);
          break;                          % continue next step
        end
        tmpuh   := uh(k); tmptk := mesh(k); % save data
        mesh(k) := min(mesh(k-1)+Q*(mesh(k)-mesh(k-1)),T_cmp);
        flag    := 1;                     % try with larger step
      end
    else
      if flag = 1                         % previous step good
        uh(k)     := tmpuh;               % recall saved data
        mesh(k)   := tmptm;
        mesh(k+1) := min(mesh(k)+(mesh(k)-mesh(k-1)),T_cmp);
        break;                            % continue next step
      else
        mesh(k) := mesh(k-1)+(mesh(k)-mesh(k-1))/Q;
        flag    := 2;                     % try with smaller step
      end
    end
  end
  if mesh(k)-mesh(k-1) < tau_min
    mesh(k)   := min(mesh(k-1)+tau_min,T);
    mesh(k+1) := min(mesh(k-1)+2*tau_min,T);
  end
end
  \end{lstlisting}
Here it is assumed that the set  $S=\{s_k\}_{k=0}^K$
 (which also includes $s_0=0$) is known a priori, while automatically finding $\{s_k\}$
 is addressed  in section~\ref{sec:heuristics} below.
  Compared to the algorithm in \cite{FrK23},
  the main changes are in the call to the residual barrier
  {\tt computeResidualBarrier}
  (where we now use the new residual barrier$\mathcal{B}(t)$) 
  and the usage of $T_{cmp}$ instead of $T$ to cut a cell at the next $s_k$
  instead of reaching over it. We highlight these changes in Algorithm~\ref{alg:mesh} by underlined text.
  To simplify the presentation, we removed the optimisations concerning the initial large factor $Q$
  given in \cite[Section 5]{FrK23} (they certainly can and should be implemented
  to get a more efficient version of Algorithm~\ref{alg:mesh}).

\begin{figure}[tb]
    \centering
       \includegraphics{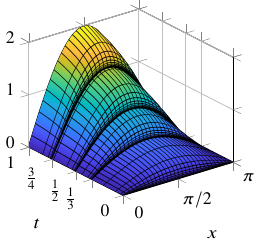}
       \includegraphics{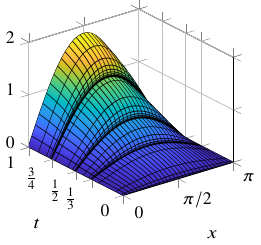}
    \caption{Numerical solution using a collocation method for \eqref{eq:ex1} (left) and \eqref{eq:ex2} (right)
             with $\alpha=0.4$, $\gamma=0.25$, $m=4$, $TOL=10^{-4}$, $Q=1.2$}\label{fig:ex12}
  \end{figure}

 \textbf{Numerical tests.} Two test examples were considered,
 both equations for $(x,t)\in(0,\pi)\times(0,1]$, subject to homogeneous initial and boundary conditions.
 For the first, we let $u$ be the solution of
  \begin{equation}\label{eq:ex1}
    (\pt_t^\alpha-\pt_x^2)u(x,t) = (H(t)+H(t-1/3)+H(t-1/2)+H(t-3/4))\sin(x).
  \end{equation}
  For the second test,  the Heaviside function components in $f$
  are smoothened to take the form $H^\gamma(t):=H(t)\, t^\gamma$,
 and we consider the solution $u$ of the equation
  \begin{equation}\label{eq:ex2}
    (\pt_t^\alpha-\pt_x^2)u(x,t) =  (H^\gamma(t)+ H^{\gamma/2}(t-1/3)+ H^{\gamma/4}(t-1/2)+ H^{\gamma/8}(t-3/4))\sin(x).
  \end{equation}
  Note, that the singularities become increasingly stronger as the
  smoothing parameter $\gamma>0$ approaches $0$ (while \eqref{eq:ex2} becomes \eqref{eq:ex1}); see Figure~\ref{fig:ex12}.

  In our numerical experiments we employed Algorithm~\ref{alg:mesh} combined with a continuous collocation method
  of order $m=4$ in time \cite{FrK23}
  and continuous cubic finite elements on 30 cells in space.
  In the residual barrier $\mathcal B$ of \eqref{eq:newBound}, we used
  all $w_k=1$, and $\lambda=0$ (which is consistent with our spatial operator).
  Figure~\ref{fig:ex12}
  shows the resulting computed solutions and the corresponding temporal meshes. We observe a strong condensing of the time steps
  immediately after each point in $S=\{0,1/3,1/2,3/4\}$, where we see local singularities  in both solutions. The first time step
  for the computed solution on the left has a width of $2.4\cdot 10^{-10}$ which is consistent with the theoretical size $TOL^{1/\alpha}\sim 10^{-10}$.
  This is also comparable with the time steps immediately after the other singularities.
  Note that the time stepping algorithm of \cite{FrK23} (with the residual barrier
  $\mathcal{R}(t)$ from \eqref {mathcal_R}) for the given parameters would
  also define a mesh and a corresponding computed solution with guaranteed error bounds, but it unnecessarily refines the mesh, and very strongly, \emph{before} the
  singularities, and, consequently,  needs twice the number of time steps in total.

  Importantly, for $TOL$ and $\alpha$ becoming smaller, the minimum time step reduces to the magnitude of $TOL^{1/\alpha}$.
  This is not a major issue near $t=0$ for the initial singularity, as here we can represent numbers
  as small as $2^{-1074}\approx 5\cdot 10^{-324}$ in double precision (non-normalised numbers). However,  with strong interior singularities,  the interior time steps may reduce in a similar way, which may lead to anothrer
  locking problem, now due to the precision being only around $2^{-53}\approx 2\cdot 10^{-16}$ (normalised numbers) for time nodes away from zero. Thus, although in theory
  the algorithm is adapting the mesh correctly, numerically one may not be able to compute the temporal mesh if $TOL^{1/\alpha}\sim 10^{-16}$.
  For the latter case,
 in the next section
 we consider an alternative  approach, which is based on a similar general idea, but requires a somewhat more intricate implementation.

  \section{Splitting and Shifting Approaches}\label{sec:Split}
  The solution in Fig.~\ref{fig:ex12}  looks like a sum of parts, each having an initial singularity at $s_k$, which is perhaps unsurprising in view of the right-hand side being split into a sum in \eqref{f_sum}.
  We can exploit this behaviour in a few ways.
 One can reduce the original problem to $K+1$ simpler problems,
  which  can then be solved using a simpler time stepping algorithm of \cite{FrK23}
  (and in parallel if the original problem is linear).
  Alternatively, one can reformulate our original equation on each of the $K$ time subintervals $(s_{k-1},s_{k})$
as an equation for $\hat t\in(0,s_k-s_{k-1})$, the latter reformulation allowing for again using the time stepping algorithm of \cite{FrK23} while avoiding the locking issues.
\smallskip

 \textbf{Splitting approach.}~%
Assuming that the original problem is linear, one can immediately split $u$ corresponding to $f$ in \eqref{f_sum} as
 $$
 u(x,t)=\sum_{k=0}^K u^k\!(x,t+s_k),
\;\;
 (\pt_t^\alpha+L) u^k\!(x,t)=f_k(x,t-s_k)\;\mbox{for~}(x,t)\in\Omega\times(0,T-s_k],
 $$
 subject to $u^0\!(\cdot,0)=u_0(x)$ and $u^k\!(\cdot,0)=0$ for $k\ge1$, and homogeneous boundary conditions.
Thus, our original problem is reduced to $K+1$ more regular problems, which can be solved, in parallel, using the adaptive algorithm from \cite{FrK23}  with a tolerance set to $w_k\cdot TOL$, which results in $K+1$ auxiliary temporal meshes and the corresponding
computed  solutions $u_{h}^k$.
The computation of the final computed solution
$ u_h(x,t)=\sum_{k=0}^K u_h^k(x,t+s_k)$,
 with a guaranteed error
  of at most $(\sum_k w_k)\cdot TOL$,
 requires the interpolation between auxiliary temporal meshes, which is a certain drawback of this, otherwise,  simple and stable approach.
\smallskip

\textbf{Shifting approach.}~%
Define the solution $u$ as a piecewise function:
$
    u_h(t)|_{[s_{k-1},s_k)}:=u_{h,k}(t+s_{k-1})
$.
  Then, for $k\in\{1,\dots,K\}$, solve the shifted problem
  \[
    (\pt_{t,-s_{k-1}}^\alpha+L)u_{h,k}(s,x)=f(t-s_{k-1},x)\;\;\text{ in }(0,s_k-s_{k-1})\times\Omega,
  \]
  with homogeneous boundary data, initial history $u_h|_{[0,s_{k-1}]}$, the shifted Caputo operator
  $\pt_{t,-s_{k-1}}^{\alpha} u(t) :=\pt_{t}^{\alpha} u(t+s_{k-1})$, and the shifted residual barrier ${\mathcal B}(t-s_{k-1})$.
  Although this idea is conceptually easier, and also applicable to non-linear operators, the fine part of the mesh
  is again shifted to zero and the mesh adaptation produces $K$ local meshes $M_k$ covering $[s_{k-1},s_k]$ each.
  But here the implementation and computation costs are higher as the problems include history terms outside the
  local mesh due to the shifted Caputo operator.

\begin{figure}[tb]
    \centering
       \includegraphics{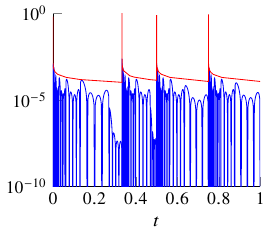}
       \includegraphics{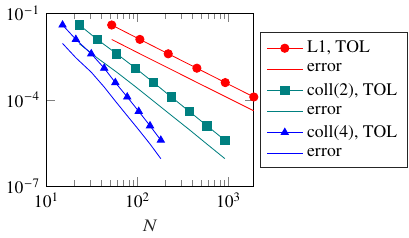}
    \caption{Example \eqref{eq:ex1} with $\alpha=0.4$, $Q=1.2$.
          Left: $L_\infty(\Omega)$ value of the residual (blue) and its barrier (red) for a collocation method with $m=4$ and $TOL=10^{-4}$.
          Right:  maximum errors vs. $TOL$.}
             \label{fig:ex1_err}
  \end{figure}

    For our numerical example \eqref{eq:ex1}, where we {again} used continuous cubic finite elements on 30 cells in space, we obtain with the second approach of shifting the time-line
    a good mesh with errors below the specified tolerances; see  Figure~\ref{fig:ex1_err} (right).
  In the left-hand picture, the mesh construction in time can be seen, where the residuals are always bounded by the barrier function $\mathcal{B}$.
  We observe quite nicely to the right of the positions of jumps in $f$ finer meshes, as the solution exhibits interior singularities there.

  \section{Finding the Problematic Positions}\label{sec:heuristics}
  Both presented approaches depend on the a-priori knowledge of where the right-hand side is non-smooth. For a fully adaptive algorithm
  it would be desirable for these positions to be found automatically. Indeed, such an automatism can be implemented in
  Algorithm~\ref{alg:mesh}. Near such a position the mesh algorithm tries to fit in increasingly smaller time steps. Thus
   catching too small time steps heuristically can be incorporated, instead of lines 35/36 of the algorithm. If such a position $s_k$ is
  recognised, it is added to the list $S$, which initially includes $s_0=0$. Then the time stepping
  should be restarted with the updated list $S$. 
  For each such restart
  one only needs to recompute the time steps starting from the last known problematic position $s_{k-1}$.
  Our experiments show that the temporal mesh produced in this way is virtually indistinguishable from a mesh with a-priori knowledge of the set $S$. But for
  the heuristics we need additional algorithm parameters, including a minimal time step leading to updating $S$ and a minimal distance to the
  last known position. In our experiments, $10^{-13}$  and $10^{-4}$, respectively, worked nicely, but these are not necessarily optimal.

  \section{Negative $\lambda$}\label{sec_lambda}
  Another reason for
 the time stepping adaptation running into locking may be  $\lambda$ becoming negative.
 Then, for example, the residual barrier from
 \eqref{mathcal_R} reduces to
  $
    \mathcal{R}(t) = -|\lambda|+\Gamma(1-\alpha)^{-1}t^{-\alpha}
  $
  and, hence, becomes negative at
  $
    s=(|\lambda|\Gamma(1-\alpha))^{-1/\alpha}.
  $
  In this very different situation, for moderate-time
  computations, exactly the same generalized residual barrier \eqref{eq:newBound}
  may still be employed, with each $s_k$
  chosen to prevent
  the $\sum_{l=0}^{k-1}$ version of $\mathcal B$ from becoming negative (or even very close to zero).
    It should also be noted that for semilinar problems, $\lambda=\lambda(t)$ is not known a priori, so the a-priori computation of such problematic points
  may be impossible; hence, Algorithm~\ref{alg:mesh} would need to be used with the heuristics of Section~\ref{sec:heuristics}.

\begin{figure}[tb]
    \centering
       \includegraphics{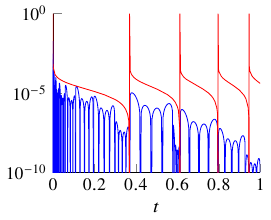}
       \includegraphics{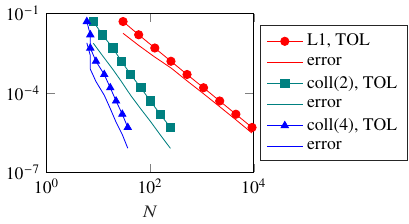}
    \caption{Example with $\lambda=-1$. Left:
             absolute value of the residual (blue) and its barrier (red) for a collocation method with $m=4$ and $TOL=10^{-4}$.
             Right: maximum errors 
             for $\alpha=0.4$, $Q=1.2$.}
             \label{fig:ex3}
  \end{figure}

  Here we test this approach, with the weights $w_k=1$ in  $\mathcal B$,
  for a very simple example without spatial derivatives:
$
    (\pt_t^\alpha-1)u(t) = f(t)$ $\text{in }(0,1)$,
  with the exact solution 
  $
    u(t)=t^{0.6}.
  $
  Figure~\ref{fig:ex3} (left)
  shows the behaviour of the residual and its barrier on the interval $[0,1]$. We observe, that the barrier $\mathcal{B}$ has five terms
  due to $\lambda=-1$ being negative. Nevertheless our algorithm is able to find the problematic positions and adjust
  the barrier function accordingly. For example, $s_1\approx 0.3695\approx\Gamma(1-\alpha)^{-1/\alpha}$. As a result we obtain a solution with a guaranteed error of
  \[
    \max_{t\in[0,1]}|u-u_h|\leq  0.3504\cdot 10^{-4}<\sum_{k} w_k\cdot TOL=5\cdot 10^{-4}.
  \]
  We  also observe in the figure that the intervals between two adjacent introductions of new bounding terms become
  smaller and smaller. This indicates, that our procedure is limited to moderate-time
  computations.
  (For arbitrarily large times, the error barrier needs to be fundamentally adjusted to allow for a positive residual barrier forall $t>0$;
  this will be addressed elsewhere in the context of more general semilinear time-fractional parabolic equations.)
 Figure~\ref{fig:ex3} (right) demonstrates that  the maximum  errors are well below the given $TOL$-values.

 \section{Computational Stability Considerations}\label{sec:stab}

\textbf{Semi-continuity of the right-hand side}\\
It turns out, that for discontinuous right-hand sides $f$ and the continuous collocation method the type of discontinuity of $f$ is important.
The continuous collocation method uses collocation points $t_k^i=t_{k-1}+c_i(t_k-t_{k-1})$, $i\in\{0,1,\dots,m\}$
where $0=c_0<c_1<\dots<c_m=1$.
Therefore, in the method (and in computing the residuals) we evaluate $f$ at $t_k^i$ for $i\in\{1,\dots,m\}$.
So it makes sense, because of $t_k^m=t_k$, to have
\[
  \lim_{t\uparrow t_k}f(t)=f(t_k)
\]
and to assume $f$ to be \textit{lower semi-continuous}. But a jump at a position $s$ in $f$ introduces a singularity in $u$ of type $(t-s)^\alpha$
and requires a mesh resolution of order $\tau\sim TOL^{1/\alpha}$ in order to have an error smaller than $TOL$.
For $\tau$ close to the precision of the computer (usually $eps=2\cdot10^{-16}$) we observe a problem in evaluating $f$ to the right of the singularity.
Here we have numerically
\[
  f(t_k^1)=f(s+c_1\tau)\stackrel{num}{=}f(s)\neq \lim_{t\downarrow s}f(t).
\]
As a consequence, the wrong value of $f$ is used, leading to an incorrect computation of either the numerical solution or the residual.
The adaptive algorithm tries to compensate for this by unnecessary refinement, which eventually locks the algorithm.

A way around this problem is to use right-hand sides $f$ that are \textit{upper semi-continuous}. Here
\[
  \lim_{t\downarrow s}f(t)=f(s)
\]
and we do not have the problem to the right of the singularity. But now we need to change the position of the last collocation point
by taking $c_m<1$. In our calculations it seems to be sufficient to choose $c_m=1-eps$.
So in a sense we introduce a mismatch in the collocation conditions. 

\textbf{Position of jumps of $f$}\\
Usually providing a-priori information about the positions of the jumps of $f$ leads to a good performance of the algorithm. 
But for $\tau$ very small we may still have problems with the last interval before the jump and the evaluation of $f$. 

A possible solution is to shift the given positions slightly to the left, e.g. by $eps$.
Now the modifications of the barrier function take effect earlier and the algorithm is more stable. Note that this repositioning
may not work if $\tau$ is much smaller than the shift.

\textbf{Evaluation of an exact solution for error calculations}\\
Even if the algorithm is stable and produces a very good mesh and approximation of the exact solution, the error computation using an exact solution may fail.
If the exact solution to the problem is known, it will contain shifted and truncated Mittag-Leffler type functions $E_\alpha$. 
These have to be evaluated at $s+\tau$ for the position $s$ of a jump and $\tau$ very small for the first cell after the jump. 
Again, numerically we evaluate the solution either at $s$, which does not yet include the jump of $f$, or at $s+eps$, 
where the exponential growth of $E_\alpha$ gives a completely wrong value.

A solution here would be to define the exact solution piecewise and to evaluate it only at local times, similar to the proposed method.
 
\bibliographystyle{plain}
\bibliography{lit}

\end{document}